\documentclass[aps,twocolumn,showpacs,groupedaddress]{revtex4}  % for review and submission
\usepackage{amssymb}
\usepackage{graphicx}
\usepackage{dcolumn}
\usepackage{bm}
\usepackage{amsmath}
\usepackage{epsfig}
\usepackage{setspace}

\begin{document} 
\title{Periodic Gabor Functions with Biorthogonal Exchange:

A Highly Accurate and Efficient Method for Signal Compression}
%\input author_list.tex       % D0 authors (remove the first 3 lines
%\affiliation{Chemical Physics Department, Weizmann Institute of Science, Rehovot, Israel}
\author{Asaf Shimshovitz and David J.~Tannor}                             % of this file 
\affiliation{Department of Chemical Physics, Weizmann Institute of Science, Rehovot, 76100 Israel}
%prior to submission, they
                             % contain a time stamp for the authorlist)
                             % (includes institutions and visitors)
\date{\today}

\begin{abstract}
We propose a new formalism for signal compression based on the Gabor basis set. By convolving the conventional Gabor functions with Dirichlet functions we obtain a periodic version of the Gabor basis set (pg). The pg basis is exact for functions that are band-limited with finite support, bypassing the Balian-Low theorem. The calculation of the pg coefficients is trivial and numerically stable, but the representation does not allow compression. However, by exchanging the pg basis with its biorthogonal basis and using the localized pg basis to calculate the coefficients, large compression factors are achieved.  We illustrate the method on three examples: a rectangular pulse, an audio signal and a benchmark example from image processing.
\end{abstract}
\maketitle
\section{Introduction}
%\subsection{Definition of the Gabor representation and why it is desirable}
In 1946, Gabor suggested using a lattice of Gaussians in time-frequency space for signal representation \cite{gabor}. (Essentially the same lattice was introduced by von Neumann in quantum mechanics in 1931 \cite{von_neumann}). The Gabor functions are constructed by shifting a synthesis function in time and frequency such that each basis function  $g_{n,m}(t)$  has the form:
\begin{eqnarray}
g_{n,m}(t)=g(t-na)\exp(jmbt) 
\label{hmn}
\end{eqnarray}
for $n,m \in Z$ where $a$ and $b$ represent time and frequency sampling intervals. By choosing the synthesis function $g(t)$ to be concentrated in both time and frequency (e.g. a Gaussian as in Gabor's original proposal), each one of the Gabor functions is localized in a different region of the time-frequency plane.  Gabor's original motivation was to use these localized $g_{n,m}(t)$ functions as a basis set for the expansion of a signal $s(t)$: 
\begin{eqnarray}
s(t)=\sum_{n,m} c_{n,m}g_{n,m}(t). 
\label{expansion}
\end{eqnarray}
For an arbitrary signal, Eq.\ref{expansion} has a solution only if $ab\leq2\pi$; for the case $ab=2\pi$, called critical sampling, the solution is unique.
This expansion is potentially very efficient for signals that contain different frequency components at different time intervals. Due to the localization of the basis set, one expects that only the basis functions in the active time-frequency regions are required for the expansion and therefore the representation of $s(t)$ should be very efficient.  

%\subsection{Statements in the literature of why it is impossible}
The Gabor basis set has been widely adopted in signal analysis, where it is known as the Short Time Fourier Transform, but it is rarely used for signal compression and reconstruction. In fact, it is generally regarded as unstable for critical sampling \cite{ingrid}. The reason for the disparity between the concept and the practice can be traced to the non-orthogonality of the basis set.  Balian and Low proved that the $g_{n,m}(t)$ cannot be made orthonormal without sacrificing the locality of $g(t)$ in either time or frequency \cite{balian,low}.  In addition, the non-orthogonality of the basis makes the calculation of the Gabor coefficients $c_{n,m}$ non trivial. In this regard, a major advance was made by Bastiaans \cite{bastians} who showed that the coefficients can be represented as the inner product between the signal $s(t)$ and a basis set $\gamma_{n,m}(t)$ that is bi-orthogonal to $g_{n,m}(t)$. However, in most cases Bastiaans's approach just transforms the non-trivial problem of finding $c_{n,m}$ to the non-trivial problem of finding $\gamma_{n,m}(t)$. For the few special cases where  $\gamma_{n,m}(t)$  can be calculated analytically, the Gabor coefficients are found to be delocalized. 

%\subsection{Discrete Gabor Theorem (wexler and raz)}
%\subsubsection{Space of periodic sequences and using implicit pbc}
About ten years after Bastiaans's work, Wexler and Raz developed a  discrete version of the Gabor expansion (henceforth `DGE') that overcomes many of the earlier difficulties \cite{wexler}. The DGE applies to sequences $s[k]$ with period $N_l$. In this formalism, the Gabor functions in Eq.\ref{hmn} turn into the discrete sequences
\begin{eqnarray}
\Tilde g_{n,m}[k]=\Tilde g[k-na]\exp(\frac{2\pi j}{N_l}mbk) 
\label{dhmn}
\end{eqnarray}
for $n=0,...,N-1$, $m=1,...,M-1$, where $a$ and $b$ are the time and frequency sampling intervals respectively and obey the relation: $Na=M b=N_l$. The synthesis sequence is now a periodic extension of some window sequence $\Tilde g[k]=\sum_l g[k+lN_l]=\Tilde g[k+N_l]$. Analogous to Eq.\ref{expansion}, the expansion of $s[k]$ is given by:
\begin{eqnarray}
s[k]=\sum_{n=0}^{N-1}\sum_{m=0}^{M-1} c_{n,m}\tilde g_{n,m}[k]. 
\label{dexpansion}
\end{eqnarray}
The coefficients $c_{n,m}$ exist for $MN\geq N_l$ and can be found by taking inner products between $s[k]$ and a 
biorthogonal basis set $\tilde\gamma_{n,m}[k]$. For critical sampling 
$NM=N_l$ the coefficients $c_{n,m}$ and $\Tilde\gamma_{n,m}[k]$ are unique and in 
contrast to the continuous-time case, can be obtained simply by solving a set of linear equations. 
%\subsubsection{Problems}
Nevertheless, some of the drawbacks of the Gabor basis persist in this formalism. To understand why, recall that the Gabor representation is expected to be efficient because of the time-frequency localization of the Gabor basis, $\tilde g_{n,m}[k]$. But this expectation is based on the assumption that only the Gabor functions in the active time-frequency window are required for the expansion Eq.(\ref{dexpansion}).
%\subsection{New modification: exchanging the roles of the Gabor and its bi-orthogonal basis}
It turns out that this key assumption is incorrect: the non-orthogonality of the basis leads to the counterintuitive phenomenon that even Gabor functions localized where the signal $s$ vanishes are required for the expansion in Eq.\ref{dexpansion}. 
%\subsubsection{As a result critical sampling is barely in use}
This delocalization of the Gabor coefficients has been overlooked in many studies \cite{daugman,zeevi,ebrahimi}; in other cases, the phenomenon was addressed by oversampling i.e. $MN>N_l$\cite{wexler,shie,bogart}. When the oversampling is 
increased sufficiently, the bi-orthogonal set $\Tilde\gamma_{n,m}[k]$ become similar to the Gabor set $\Tilde g_{n,m}[k]$ and the coefficients become more localized. However, oversampling is inefficient in the sense that a time-frequency region that was previously overlapped by one Gabor function in is now overlapped by many, and all of them need to be taken into account in the expansion in Eq.\ref{dexpansion}. 

%\subsection{Difficulties in Gabor leads to different approach-wavelet}
As a result of these difficulties, a different approach for localized time-frequency expansion was developed called the ``wavelet'' expansion \cite{ingrid}, in which the time-frequency localization is achieved by scaling the synthesis function in addition to shifting. The great advantage of this approach is that it allows one to construct a localized basis set that is orthogonal, and therefore bypass the problems of the Gabor expansion.

%\subsection{New method will be developed in this paper - \underline{\textit{pg(b?)}}}
In this work, we return to the original Gabor formulation and introduce two modifications that make the representation both stable and highly efficient.  The method can also be combined with wavelet or pyramidal scaling, although this will be addressed in a separate publication.  The first modification is that that Gabor functions are `convolved' with a Dirichlet kernel (a periodic sinc function), giving a periodic, Gabor-like basis that is guaranteed to be stable since it spans exactly the time-frequency space to which Nyquist's theorem applies. (In the final expressions, all that appears are the original Gabor functions at the sample points.) Second, the non-locality of the coefficients is overcome by exchanging the role of the Gabor basis and its bi-orthogonal basis set. Specifically, the localized periodic-Gabor functions are used to calculate the coefficients rather that as the basis. We call the combination of the two innovations the ``periodic Gabor method with biorthogonal exchange", or pgb. In \cite{asaf,norio} we applied the pgb to quantum mechanics with great success \footnote{In the context of quantum mechanics we call the method pvb for periodic von Neumann basis with biorthogonal exchange.}; here we show here that the method can provide large savings for signal compression as well.   
\section{The Periodic Gabor Basis}
%\subsection{Gabor}
The original proposal by Gabor was to use  
Gaussian functions on a time-frequency grid with critical sampling \cite{von_neumann,gabor}. Under these conditions Eq.\ref{hmn} turns into:
\begin{eqnarray}
g_{n,m}(t)=\left(\frac{2\alpha}{\pi}\right)^\frac{1}{4}\exp\left(-\alpha(t-t_n)^{2}+j\omega_m(t-t_n)\right)
\label{G}
\end{eqnarray}%
where $n$ and $m$ are integers. Each basis function is a Gaussian centered
at $(t_n,\omega_m)=(na+t_0,\frac{2\pi m}{a}+\omega_0)$ in time-frequency space, where $t_0$ and $\omega_0$ are arbitrary shifts.  
The parameter
$\alpha=\frac{\sigma_\omega}{2\sigma_t}$ controls the FWHM of each Gaussian in $t$ and
$\omega$ space. Taking $\Delta t =a$ and $\Delta \omega
=2\pi/a$ as
the spacing between neighboring Gaussians in $t$ and $\omega$ space
respectively, we note that $\Delta t \Delta\omega=2\pi$
so we have exactly one basis function per unit cell in
phase space, so-called critical sampling. Assuming that the grid extends over the infinite time-frequency grid, the Gabor basis is complete as shown  in \cite{perelomov}.

However, in any numerical calculation, $n$ and $m$ take on only a finite number of values. Consider a basis consisting of $N$ Gaussian functions $\lbrace g_i(t)\rbrace$, $i=1...N$, one per unit cell. Since the size of one Gabor unit cell is $2\pi$, the phase space area covered by the truncated Gabor lattice is:   
\begin{eqnarray}
S^{\rm Gabor}=2\pi N.
\label{AREA}
\end{eqnarray}
\begin{figure}[h]
\begin{center}
\includegraphics [width=7cm]{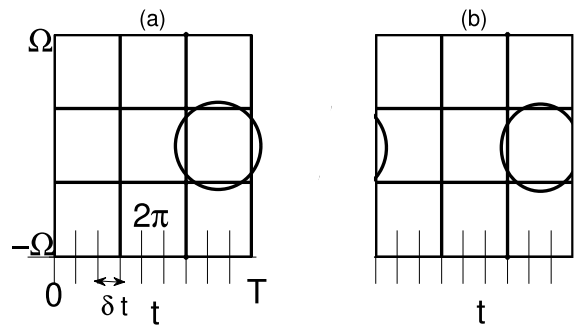}
\end{center}
\begin{center}
\caption{\footnotesize{(a) $N=9$ coordinate grid points and $N=9$ Gabor unit cells cover
the same area in phase space, $S=2\pi N$. Shown also is a typical Gabor function.  Note that its boundary conditions are not appropriate for the rectangular area. (b) The periodic Gabor (pg) basis is a complete set  for the truncated space. The pg basis
functions are loosely speaking, periodic Gaussians whose centers are located at the center of each unit
cell. }}
\label{VN}
\end{center}
\end{figure}
However, it turns out that there is a difference between the truncated Gabor lattice and the truncated Gabor \textit{basis}: The latter is not complete for the Hilbert space defined by the rectangular area $S^{\rm Gabor}=2\pi N$ because the boundary conditions are not appropriate.  

%\subsection{Dirichlet kernel}
In order to develop a truncated Gabor-like basis set that is complete on the rectangular area in Eq.\ref{AREA} (Fig.\ref{VN}) we introduce a variant on the Gabor basis set with periodic boundary conditions. Consider $x(t)$ to be  a periodic signal with period $T$ and band limited in $\Omega$. By sampling $x(t)$ at the Nyquist rate   $\delta t=\frac{\pi}{\Omega}$ at $N=\frac{T}{\delta t}$ sampling points $\{t_i\}$, $x(t)$ can be expressed in terms of its samples by:
\begin{eqnarray}
x(t)=\sum_{i=1}^{N}x(t_i)D(t-t_i) 
\label{sample_rec}
\end{eqnarray}
where $D(t)$ are the Dirichlet or periodic sinc functions\cite{nandi}:
\begin{eqnarray}
D(t)&=&\frac{\sin(N\Omega t/2)}{N\sin(\Omega t/2)}
\nonumber \\
&=&\frac{1}{N}\sum_{m=-(N-1)/2}^{(N-1)/2}\exp(jm\Omega t).
\end{eqnarray}
Therefore the set $\lbrace D_i(t)\rbrace=D(t-t_i)$ for $i=1,...,N$   comprises a complete basis set for the time-frequency rectangular space $S^{\rm{D}}=2T\Omega=2N\delta t\frac{\pi}{\delta t}=2\pi N$. 

%\subsection{pg}
By combining the Dirichlet and Gabor basis functions we
can generate a finite ``Gabor-like'' basis set that is complete on the
truncated space given in Eq.\ref{AREA}: 
\begin{eqnarray}
\Tilde g_m(t)=\sum_{i =  1}^{N}g_m(t_i)D(t-t_i)
\label{dvn}
\end{eqnarray}
for $m=1,...,N$. The new basis set is in some sense, the Gabor functions with periodic boundary conditions (henceforth, the periodic Gabor or `pg' basis). We can write eq.(\ref{dvn}) in matrix notation as: $\mathbf{\Tilde {G}}=\mathbf{D}\mathbf{G}$ where $G_{ij}=g_j(x_i)$. By taking the width parameter $\alpha=\frac{\Delta\omega}{2\Delta t}$ (and locating a Gaussian function at the center of each unit cell) we
can guarantee that the Gabor functions have no linear dependence and
that the matrix $\mathbf{G}$ is invertible, that is $\mathbf{\Tilde G G^{-1}=D}$.

The invertibility of $\mathbf{G}$ implies that both the Dirichlet and the pg bases span the same space (Fig.\ref{VN}).  
Therefore, $x(t)$ can be expanded as:
\begin{eqnarray}
x(t)=\sum_{m =  1}^{N}\Tilde g_m(t) c_m.
\label{psi}
\end{eqnarray}
To find the coefficients $a_m$ we first define the overlap matrix, $\mathbf{S}$, as the inner product between pg functions: 
\begin{eqnarray}
S_{ij}&=& \langle\Tilde g_i , \Tilde
g_j\rangle=\int_{0}^{T_0} \Tilde g_i^{*}(t)\Tilde g_j(t)dt
\nonumber \\
&=&\sum_{n =  1}^{N}\sum_{m =  1}^{N} g_i^{*}(t_n)
g_j(t_m)\int_{0}^{T_0} D^{*}(t-t_n)D(t-t_m)dt
\nonumber \\
&=&\sum_{n =  0}^{N-1} g_i^{*}(t_n)g_j(t_n)
\label{ip}
\end{eqnarray}
or
\begin{equation}
\mathbf{S=G^{\dag}G}.
\label{S}
\end{equation}
Using the completeness relationship for non-orthogonal bases \cite{tannor_book3}, 
$x(t)$ can be expressed as
\begin{equation}
x(t)=\sum_{n =  1}^{N} \sum_{m =  1}^{N}
\Tilde g_m(t)(S^{-1})_{mn} \langle \Tilde
g_n, x\rangle,
\label{psi2}
\end{equation}
where, similar to Eq.\ref{ip}, the inner product 
\begin{equation}
\langle \Tilde
g_n, x\rangle=\sum_{n =  1}^{N} g_i^{*}(t_n)x(t_n).
\label{ip_g_signal}
\end{equation}
Comparing Eq. (\ref{psi2}) with Eq.(\ref{psi}) we find that 
\begin{equation}
c_m=\sum_{n =  1}^{N} (S^{-1})_{mn}\langle \Tilde
g_n, x\rangle \label{pvn_coef}
\end{equation}
for $m=1,...,N$. Although the $\mathbf{S}$ matrix is localized, $\mathbf{S}^{-1}$ is not; hence, Eq. \ref{pvn_coef} indicates that basis functions $\Tilde g_m$ that are distant in time-frequency from the signal $x$, still contribute to the expansion because of $\mathbf{S}^{-1}$. 

Equation \ref{pvn_coef} can be written in an alternative way:
\begin{equation}
x(t)=\sum_{m =  1}^{N} \Tilde g_m(t) \langle \Tilde b_m, x\rangle,
\label{psi2b}
\end{equation}
where
\begin{eqnarray}
\langle \tilde b_m, \equiv \sum_{n =  1}^{N}(S^{-1})_{mn}
\langle \tilde g_n,
\label{fvn}
\end{eqnarray}
for $m=1,...,N$. We will call the $\{ \tilde b_m \}$ basis the biorthogonal Gabor or `bg' basis since the pg and bg bases taken by themselves are 
non-orthogonal but are orthogonal to each other.  
This can be shown easily by:
\begin{eqnarray}
\langle \Tilde b_i, \Tilde g_j \rangle&=&\sum_{n =  1}^{N} b^{*}_i(t_n)g_j(t_n)
\nonumber\\
&=&\sum_{m =  1}^{N} \sum_{n =  1}^{N} g^{*}_m(t_n)g_j(t_n)(S^{-1})_{im}
\nonumber\\
&=&\sum_{m =  1}^{N} S_{mj}(S^{-1})_{im}=\delta_{ij}.
\end{eqnarray}
%Eq.13-14 provides great advantages over the G basis set
This biorthogonality property was used to great advantage by Bastiaans on the infinite basis set, but here the biorthogonality property holds on a truncated rectangular region in the time-frequency space.

So far the pg basis set provides two great advantages over the traditional Gabor basis set. First, in the original Gabor basis $N$ functions give a very inaccurate expansion of $x(t)$\cite{low}, whereas the periodicity of the pg basis makes the expansion exact for functions that are band-limited with finite support, bypassing the Balian-Low theorem.  Second, the coefficients in the original Gabor basis require knowledge of the bi-orthogonal basis, which in most cases is non-trivial to calculate. In the pg basis, no knowledge of the bi-orthogonal basis is necessary (Eq.\ref{pvn_coef}). Nevertheless, the compression is still problematic since, as discussed at Eq. \ref{pvn_coef}, basis functions that are distant in time-frequency from the signal still contribute to the expansion. 
\section{Achieving Compression: The Periodic Gabor Basis with Biorthogonal Exchange}
%\subsection{Bi-orthogonal basis set}
We now turn to the issue of compression.  Gabor's original idea was to exploit the time-frequency localization of the basis functions to reduce the size of basis.  This intuitively appealing ideas translates to the statement that a significant fraction of the basis will fulfill the relation: $\langle\Tilde g_n,x\rangle=0$, $n=1,...,M$. However, there is an important subtlety --- we cannot simply eliminate the states $\Tilde{g}_n$, since as noted at Eq. \ref{pvn_coef}, even pg functions that are remote from the signal can contribute to the expansion. Nevertheless, we can take advantage of the near-vanishing $\langle \Tilde g_n, x\rangle$, if we exchange the roles of the bg and pg bases, allowing the delocalized bg functions to serve as the basis set for $x(t)$. Then Eq.\ref{psi2} becomes: 
\begin{eqnarray}
x(t)=\sum_{n =  1}^{N}  b_n(t) d_n=\sum_{n =  1}^{N}  b_n(t) \langle \Tilde g_n, x\rangle. 
\label{bvn_rec}
\end{eqnarray}
%where $c_n$ is given by Eq.\ref{ip_g_signal}.
By assumption, $M$ of the $\langle \Tilde g_n, x\rangle$ are zero, hence in order to
represent $x(t)$ in the bg basis set  we need only $N'=N-M$
basis functions. To emphasize, although the basis $\{b_n (t)\}$ has played a central role in Gabor theory since 1980, it has always been used for finding the  \textit{coefficients} of the Gabor basis.  Equation \ref{bvn_rec} is a radical departure from the existing literature since the roles of the (periodic) Gabor and the biorthogonal (periodic)Gabor bases have been exchanged, with the bi-orthogonal functions now actually serving as the \textit{basis} for the expansion of the signal and the localized original functions used for the calculation of the coefficients.  Hence the name of the method,``periodic Gabor theory with biorthogonal exchange", or pgb.
   
%\subsection{Porat's correction}
As discussed above, a large number of the coefficients in the pgb method are expected to be close to zero. After eliminating the small coefficients by quantization or entropy compression, one can use the remaining coefficients and the corresponding bg functions to obtain an approximate reconstruction of the original signal. We have found that the accuracy of the reconstruction, although quite good, can be further improved using a method developed by Porat\cite{porat}. The first step in Porat's method is to find a basis that is biorthogonal to the contracted bg set, i.e. by repeating the procedure of Eq.\ref{ip} and Eq.\ref{fvn} but with the role of the $N$ pg functions played by the $N'$ bg set. In particular, the $\mathbf{S}$ matrix is now defined in terms of the contracted set only. Then the desired coefficients are obtained by the inner product between this new biorthogonal set and the original signal. We will show results using the Porat correction below.
\section{Discrete formulation of the pgb method}
The pgb method can be formulated for discrete sequences as well as for continuous signals. Consider a discrete sequence $x[n]$ of length $N$, which can be viewed as samples taken at $t_i= 0,1,...,N-1$ of an $N$ periodic signal $x(t)$ that is band limited in the interval $[-\pi,\pi]$. Let the Gabor grid consist of $N_t$ and $N_{\omega}$ Gaussians in $t$ and $\omega$ spaces respectively such that $N_tN_{\omega}=N$. The parameters that define the Gabor functions in Eq.\ref{G} are $\alpha=\frac{\Delta\omega}{2\Delta t}=\frac{\pi}{(N_\omega)^2}$, $t_0=-\frac{\Delta t}{2}$, $\omega_0=-\frac{2\pi+\Delta\omega}{2}$ , $n=1,...,N_t$ and $l=1,...,N_\omega$. The pgb representation of the sequence is given by: 
\begin{eqnarray}
x(n)=\sum_{i =  1}^{N}  b_i(n) d_i, 
\label{dbvn}
\end{eqnarray}
or in matrix notation:
\begin{equation}
\rm{x}\mathbf{=B}\rm{d},
\end{equation}
where x is the length $N$ column vector of the sequence of values, $\mathbf{B}$ is the $N \times N$ matrix with elements $B_{ni}=b_i(n)$ and d is the coefficient vector given by Eq.\ref{ip_g_signal}. Using Eq.\ref{fvn} and the fact that $\Tilde g(n)=g(n)$, we find that:
\begin{eqnarray}
b_i(n)=\sum_{j =  1}^{N} g_j(n)(S^{-1})_{ji}. 
\end{eqnarray}
or in matrix form:
\begin{equation}
\mathbf{B=GS^{-1}=(G^{\dagger})^{-1}},
\label{bgdagger}
\end{equation}
where we have used Eq. (\ref{S}). 
Equation \ref{bgdagger} can be written succinctly as $\mathbf{BG^{\dagger} =1}$. The related equation $\mathbf{GB^{\dagger}=1}$ is also correct, but has a completely different physical significance --- it corresponds to the expansion of the signal in terms of the original Gabor basis, without the exchange of roles with the biorthogonal basis; in other words it corresponds to a discrete pg representation.  It turns out that this discrete pg representation is fully equivalent to the DGE. In fact, the equation $\mathbf{GB^{\dagger}=1}$ appears in \cite{wexler}, albeit for just one column of the matrix.
\section{Results}
\begin{figure}[h] 
\begin{center}
\includegraphics [width=8cm]{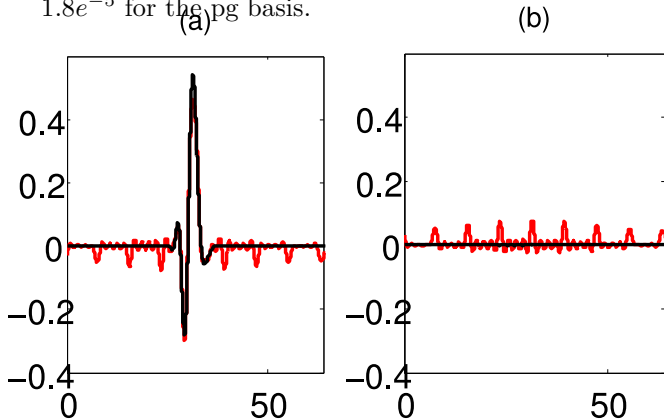}
\end{center}
\begin{center}
\caption{\footnotesize{(a). Reconstruction of a chirped Gaussian function using 64 standard Gabor(red) and 64 pg (black) functions. The pg result is virtually indistinguishable from the exact result. (b) Deviation between the reconstructed and the original signal using standard Gabor (red) and pg (black) functions.}}
\label{continous_rec}
\end{center}
%\vspace{.cm}
\end{figure}

%\subsection{pg vs standard Gabor}
We first address the difference in accuracy between the standard Gabor basis and the pg basis set without exchange. %As discussed above, the pg basis when applied to a discrete signal is equivalent to the DGE.  
Consider the chirped Gaussian $s(t)=(\frac{2\rm{Re}(\alpha)}{\pi})^{\frac{1}{4}}\exp(-\alpha(t-t_0)^2+iw_0(t-t_0))$ with $\alpha=0.135(1+i)$, $t_0=31.5$ and $w_0=\frac{\pi}{3}$. We will compare the reconstruction using 64 Gabor and 64 pg functions. Both basis sets consist of Gaussian functions on an $8 \times 8$ time-frequency lattice. The Gabor coefficients were calculated using the procedure in \cite{porat}. Figure \ref{continous_rec}(a) shows the real part of the original signal and the reconstructions obtained using both bases. Figure \ref{continous_rec}(b) shows the difference between the reconstructed and the original signal. The norm of the error is $0.0347$ for the Gabor basis and $1.8e^{-5}$ for the pg basis. 
%\begin{figure}[h] 
%\begin{center}
%\includegraphics [width=8cm]{continus_chirp_gaussian_reconstruction_16_functions_pvn_gabor.ps}
%\end{center}
%\begin{center}
%\caption{\footnotesize{(a). Reconstruction of a chirped Gaussian function using 64 standard Gabor(red) and 64 pg (black) functions. The pg result is virtually indistinguishable from the exact result. (b) Deviation between the reconstructed and the original signal using standard Gabor (red) and pg (black) functions.}}
%\label{continous_rec}
%\end{center}
%\vspace{.cm}
%\end{figure}

%\subsection{The pgb basis vs discrete Gabor expansion}
We now turn to compression using the pgb method. (By compression we mean using only a small subset of the coefficients for the reconstruction.) We present three examples.  The first is a 1-d rectangular pulse, an important exmaple for signal processing. The second is the 2-d ``Lena" picture, a benchmark for image processing.  The third is an audio signal corresponding to a `splat' sound.  For each example, we compress the signal using a small subset of the coefficients and present the reconstruction from this subset. 

\begin{figure}[h]
\begin{center}
\includegraphics [width=9.0cm,height=4.7cm]{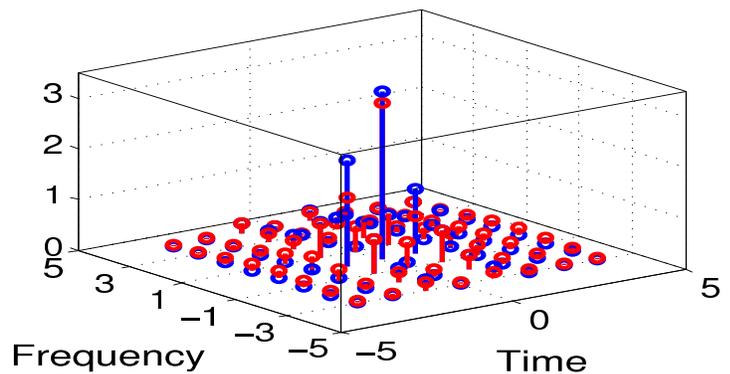}
\end{center}
\begin{center}
\caption{\footnotesize{The coefficients for representing a rectangular pulse in the DGE (red) and in the pgb method (blue).}}
\label{figure:spectrum}
\end{center}
\vspace{-.5cm}
\end{figure}

%\subsubsection{Rectangular pulse}
The rectangular pulse was studied previously by Wexler and Raz\cite{wexler}. Following those authors we consider a rectangular pulse represented by a discrete sequence of 64 points, compressed using 25 coefficients. 
\begin{figure}[h] 
\begin{center}
\includegraphics [width=5.3cm,height=3.5cm]{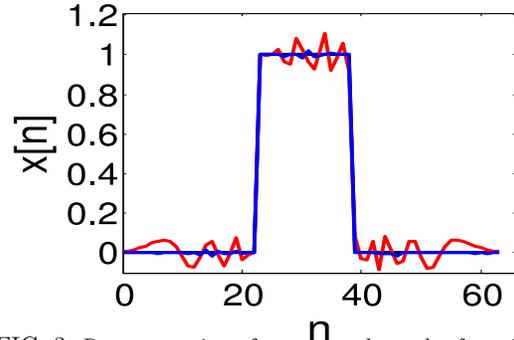}
\end{center}
\begin{center}
\caption{\footnotesize{Reconstruction of a rectangular pulse from 25 coefficients using the DGE method (red) and the pgb method (blue).}}
\label{pulse}
\end{center}
%\vspace{.cm}
\end{figure}
The reconstruction obtained with the pg basis is essentially identical to the one reported by Wexler and Raz using the DGE \cite{wexler} (see the discussion below Eq. \ref{bgdagger}).  However, with the pgb method the reconstruction is improved by essentially an order of magnitude ($L^2$ norm of the error is 0.3989 in the DGE and 0.0469 in the pgb method).
Comparison of the spectrum in the DGE and the pgb methods (Fig.\ref{figure:spectrum}) demonstrates the much higher localization of the pgb coefficients, specifically in the time coordinate (the signal is inherently delocalized in frequency).
%\subsubsection{Lena Image}
%\begin{figure}[h]
%\begin{center}
%\includegraphics [width=3.9cm]{lena.ps}
%\end{center}
%\begin{center}
%\caption{\footnotesize{Original lena image contains $512\times 512$ pixels.}}
%\label{lena}
%\end{center}
%\vspace{-.5cm}
%\end{figure}

\begin{figure}[h]
\begin{center}
\includegraphics [width=9.0cm]{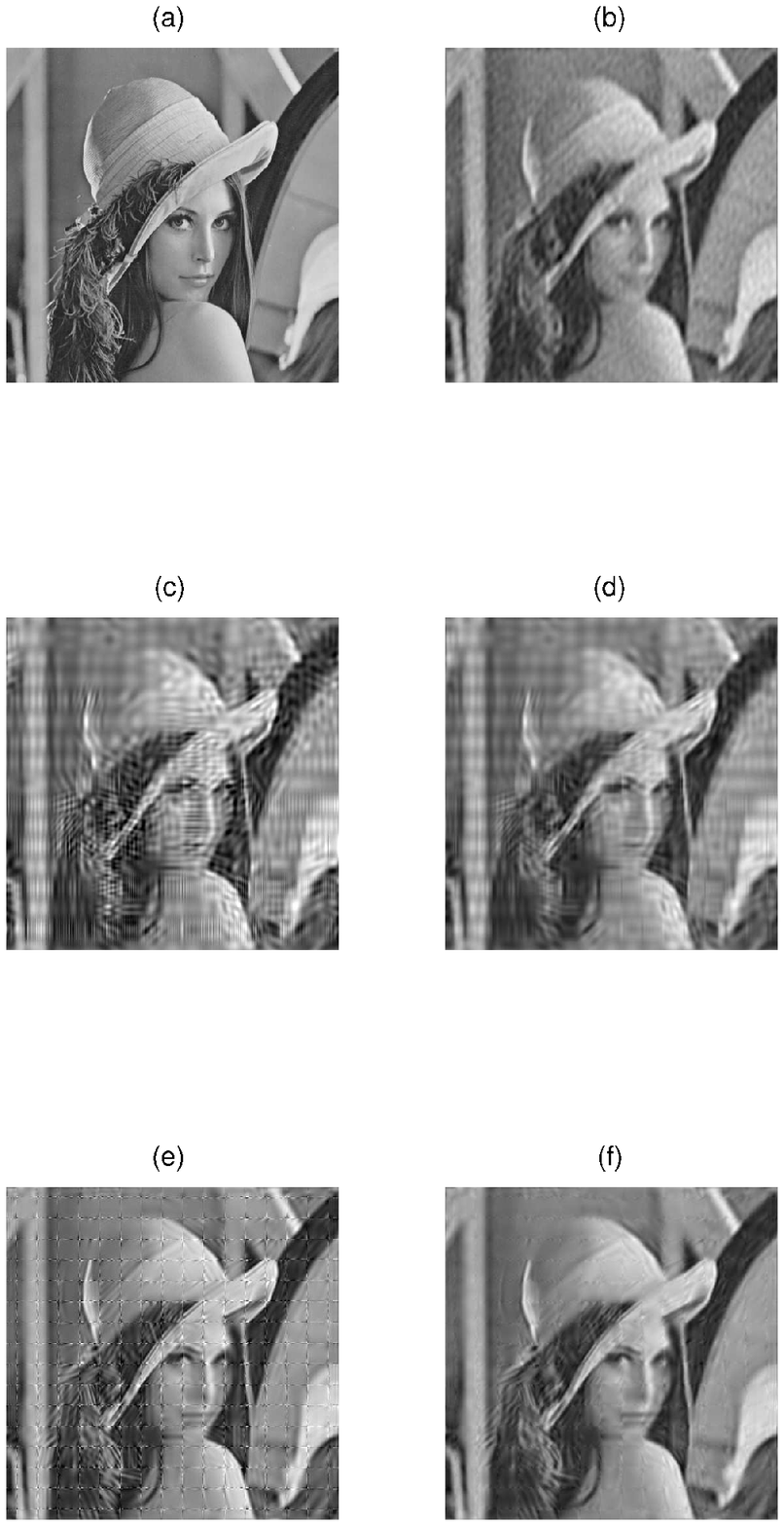}
\end{center}
\begin{center}
\caption{\footnotesize{Reconstruction of the Lena image using about $\approx$ 1\% (2621) of the coefficients. (a) original picture. (b) DFT. (c) DGE method. (d) DGE with Porat's correction. (e) pgb method. (f) pgb method with Porat's correction. }}
\label{lena_rec}
\end{center}
\vspace{-.5cm}
\end{figure}

Several works have extended the DGE to two dimension in order to apply it to image compression\cite{ebrahimi,Li,andras,si}. Although the 2-d DGE reduces the entropy of an image, the entropy can be reduced much further using the pgb method. We will demonstrate this on
the Lena image (Fig. \ref{lena_rec}(a)), which is a benchmark for image compression. The original image contain $512\times 512$ pixels. Figure \ref{lena_rec} shows the reconstruction of the image using only the 2621 largest coefficients ($\approx$ 1$\%$ of the coefficients), using the DFT (b), the DGE method (c) and the pgb method (e). The mean-square-error of the reconstructed image is 417 in the DGE and 327 in the pgb method, but the difference in visual quality is far more significant than these numbers indicate.
Figure \ref{lena_rec}(d,f) shows the reconstruction of the image using 2621 coefficients when Porat's correction is applied to both methods. The mean-square-error of the reconstructed image is 279 in the DGE and 144 in the pgb method.  Again, the difference in visual quality is far more significant than these numbers indicate. The pgb result with Porat's correction is better than the DFT (mean-square-error of 173) with the same number of coefficients (Fig. \ref{lena_rec}(b)).
%\begin{figure}[h]
%\begin{center}
%\includegraphics [width=9.0cm]{99_compression_periodic_bvn_vs_periodic_pvn_porat.ps}
%\end{center}
%\begin{center}
%\caption{\footnotesize{Reconstruction of the Lena image using 2621 coefficients with Porat's correction for (a) the standard Gabor method (b) the pgb method.}}
%\label{lena_rec_porat}
%\end{center}
%\vspace{-.5cm}
%\end{figure}
%\begin{figure}[h]
%\begin{center}
%\includegraphics [width=3.5cm]{99_compression_fourier.ps}
%\end{center}
%\begin{center}
%\caption{\footnotesize{Reconstruction of the Lena image using 2621 coefficients with Porat's correction for (a) the standard Gabor method (b) the pgb method.}}
%\label{lena_rec_porat}
%\end{center}
%\vspace{-.5cm}
%\end{figure}
%\subsection{pgb vs Time-Frequency Representaion}

\begin{figure}[h]
\begin{center}
\includegraphics [width=5cm]{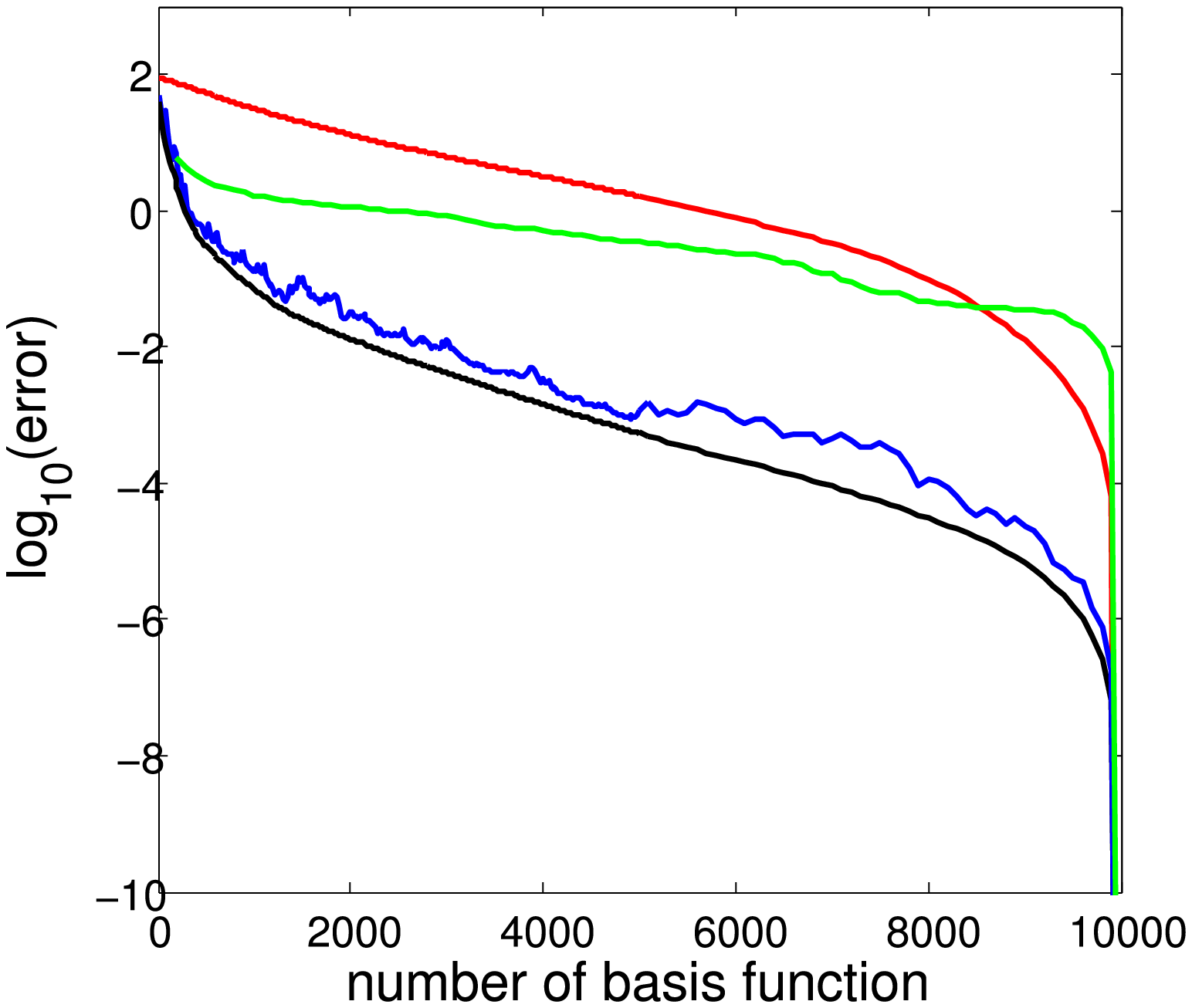}
\end{center}
%\vspace{.6cm}
\begin{center}
\caption{\footnotesize{The norm of the error of the reconstructed signal as a function of the number of basis functions using the DFT (red), the DGE with Porat's correction (green), the pgb (blue), and the pgb with Porat's correction (black).}}
\label{comp_splat}
\end{center}
%\vspace{-1.2cm}
\end{figure}

%\subsubsection{Compression of an audio signal}
Finally, we consider an example of an audio signal, which demonstrates the advantage of the pgb representation over the time and frequency representations.
%A more realistic example is given by the bvN representation of a simple audio signal. 
The original audio signal is sampled at 10000 sampling points. We used the pgb method on a grid of 100*100 unit cells in time-frequency space.
\begin{figure}[h]
\begin{center}
\includegraphics [width=6cm]{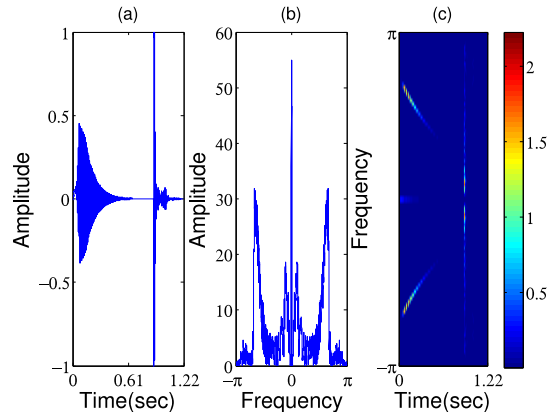}
\end{center}
\begin{center}
\caption{\footnotesize{The splat signal in time, frequency and pgb representations. }}
\label{splat_representation}
\end{center}
\end{figure}   
Figure(\ref{splat_representation}) shows the signal in the time, frequency and pgb representations. We reconstructed the signal for various basis set sizes and calculated the norm of the error (the difference between the original and the reconstructed signal) using both the original pgb coefficients also the pgb coefficients after Porat's correction. The results are shown in Fig. (\ref{comp_splat}). For comparison, we show results retaining the same number of DFT coefficients. 
Clearly, the pgb method allows much more compression than the DFT, but introduces some residual roughness.  This roughness may be removed using Porat's correction, although the latter involves some additional computational cost.
\section{Summary}
We introduced a new basis set called the periodic Gabor or pg basis which is effectively a `convolution' of the Gabor and Dirichlet (periodic sinc) functions. The original Gabor basis is complete for the infinite space but incomplete, and actually unstable for a truncated rectangular region in time-frequency space. In contrast, the pg basis provides a complete and stable representation for a truncated rectangular area in time-frequency space, i.e. for the Hilbert space spanned by functions that are band-limited with finite support. Moreover, the coefficients in the pg representation may be calculated by simple matrix multiplication (Eq.\ref{pvn_coef}). 

We noted that although the pg functions are localized, the corresponding biorthogonal bg basis functions that determine the coefficients are not.  By exchanging the roles of the pg and bg bases, (``periodic Gabor method with bi-orthogonal exchange" or `pgb') the coefficients become localized, leading to large compression factors. The pgb formalism is trivially adapted to the case of discrete signals. We provided three examples of finite sequences, and showed the large compression factors achievable using the pgb coefficients. 

Although we have not exploited it here, Eq. (\ref{dvn}) combined with Eq. (\ref{AREA}) allows the freedom to choose Gaussian basis functions of different widths and spacings. That means that in principle we can combine the time-frequency shifting with a scaling transformation, similar to what is done in the wavelet approach and in the pyramidal Gabor function. In fact, the formalism is not restricted to Gaussians --- it can be extended to any localized basis functions as long as the unit cell area obeys $\Delta t\Delta\omega=2\pi$. More research will be required on how to optimize the basis for a given signal.  

We fully expect that the pgb method and its generalizations, with its combination of simplicity, flexibility, accuracy, stability and efficiency, will provide a competitive alternative to wavelet methods.

This work was supported by the Israel Science Foundation Grant 807/08 and the US-Israel Binational Science Foundation. It was made possible in part by the historic generosity of the Harold Perlman family.

\end{document}